\documentclass[english]{article}
\usepackage[T1]{fontenc}
\usepackage[latin9]{inputenc}
\usepackage[letterpaper]{geometry}
\usepackage[active]{srcltx}
\usepackage{amsthm}
\usepackage{amsmath}
\usepackage{amssymb}
\PassOptionsToPackage{normalem}{ulem}
\usepackage{ulem}

\usepackage[all]{xy}
\usepackage{pb-diagram}
\makeatletter
\numberwithin{equation}{section}
\numberwithin{figure}{section}
\theoremstyle{plain}
\newtheorem{theorem}{\protect\theoremname}[section]
  \theoremstyle{plain}
  \newtheorem{proposition}[theorem]{\protect\propositionname}
  \theoremstyle{plain}
  
\newtheorem{remark}[theorem]{Remark}
\newtheorem{defi}[theorem]{Definition}
\newtheorem{example}[theorem]{Example}
\@ifundefined{date}{}{\date{}}
\usepackage{tikz-cd}

\makeatother

\usepackage{babel}
  \providecommand{\corollaryname}{Corollary}
  \providecommand{\propositionname}{Proposition}
\providecommand{\theoremname}{Theorem}

\def\R{\mathbb{R}}
\def\C{\mathbb{C}}
\def\H{\mathbb{H}}
\def\K{\mathbb{K}}
\def\a{\overline{a}}

\def\az{\overline{a_0}}
\def\aa{\overline{a_1}}

\def\bz{\overline{b_0}}
\def\ba{\overline{b_1}}
\def\bb{\overline{b_2}}

\begin{document}
	\title{Nonsingular Bilinear Maps Revisited}
	\author{Carlos Dom\'{i}nguez \footnote{Academ\'ia de Matem\'aticas, Unidad Interdiciplinaria de Ingenier\'{i}a Campus Guanajuato, Instit\'uto Polit\'ecnico Nacional,  Av. Mineral de Valenciana 200, Fracc. Industrial Puerto Interior, C. P. 36275, Silao de la Victoria, Guanajuato, México. \newline Email: cdomingueza@ipn.mx } , Kee Yuen Lam  \footnote{Department of Mathematics, University of British Columbia, Vancouver, B.C.,V6T 1Z2, Canada. \newline Email: lam@math.ubc.ca }}
	
	\maketitle

	\section{Introduction} \label{intro}
	
	In this paper bilinearity always means $\R-$bilinearity. A bilinear map $\varPhi:\R^r \times \R^s\to \R^n$ is nonsingular if $\varPhi(\overrightarrow{a},\overrightarrow{b})=\overrightarrow{0}$ implies $\overrightarrow{a}=\overrightarrow{0}$ or $\overrightarrow{b}=\overrightarrow{0}.$ These maps generalize the multiplication of the classical division algebras $\R,\C,\H $ and $\K$ of real, complex, quaternion or octonion numbers, which correspond to the cases $r=s=n=1,2,4$ and $8$ respectively. The study of such maps by topological methods began with Hopf \cite{Hopf} and Stiefel \cite{Stiefel}, leading eventually to their applications to embedding and immersion of real projective spaces $\R P^{r-1}$ into Euclidean Space \cite{Ginsburg,Milgram,Steer}. Applications to homotopy groups of spheres can be found in \cite{Kee3}.

    The problem, for what triples $(r,s,n)$ can there exist a nonsingular bilinear map $\varPhi,$ remains unsolved up to the present day. See \cite{Kee5} for an overall discussion. From past experience non-existence results will have to involve increasingly sophisticated tools in algebraic topology \cite{Kitchloo-Wilson}, while existence results can be obtained through skillful constructions via algebra. In particular, use of the octonions $\K,$ with its ``restricted associativity property'' listed in Section 2, has led to constructions by Lam \cite{Kee}  and Adem \cite{Adem,Adem1,Adem2}.
    
    More than 47 years have passed since \cite{Adem2} without any new example of nonsingular bilinear maps appearing in literature. In particular, an open question about commutators in $\K,$ posed in the introduction of \cite{Adem2}, remains unanswered. The purpose of this paper is to construct a new family of nonsingular bilinear maps, in Section 6 below, and to comment on their topological implications. In particular, we answer Adem's open question in the affirmative.

    The second author would like to dedicate this work to the memory of professor Elmer Rees (1941-2019), a highly esteemed colleague and long time friend. Elmer's paper \cite{Rees}, in particular, has been a source of inspiration for the present article. Both authors are indebted to the late professor Jos\'e Adem, whose papers \cite{Adem,Adem1,Adem2} set the stage for Theorem \ref{Thmmain} bellow.

\section{Restricted associativity of the octonions $\K$}\label{AsoOct}

For elements $a,b,c, \cdots$ in  $\K$ we shall use $\R(a,b,c,\cdots)$ to denote the subalgebra they generate. The following properties are well-known, and shall be frequently used in the sequel. 

\begin{enumerate}
	\item $\R(a)$ always contains the conjugate $\overline{a}.$
	\item If $ab=ba,$ then $\R(a,b)$ is commutative as well as associative. There exists then $d\in \K$ such that $\R(a,b)=\R(d)=$ a field, isomorphic either to $\R$ or to $\C.$ 
	\item (restricted associativity). If $ab=ba,$ then for any $c\in \K, \R(a,b,c)=\R(d,c)$ is associative, and. Restricted associative laws hold:
	
	$$a(bc)=(ab)c; \;\;(ca)b=c(ab).$$
	
	\item In particular, since $a$ always commutes with $\overline{a},$ one has $$a(\a c)=(a\overline{a})c=|a|^2c=(ca)\overline{a}.$$   
\end{enumerate}

\section{Twisted polynomial multiplication and bilinear maps}\label{twistedpolmult}

A primary source of nonsingular bilinear maps is polynomial multiplication. In this paper we shall mainly deal with the polynomial ring $\Lambda[X]$ with coefficient ring $\Lambda=\K,$  de-emphasizing the cases $\Lambda= \H,\C$ or $\R.$ If

      $$p(X)=a_0+a_1X+\cdots + a_rX^r$$
      $$q(X)=b_0+b_1X+\cdots +b_sX^s,$$
      then $p(X)q(X)$ is traditionally defined to be 
      $$p(X)q(X)=c_0+c_1X+c_2X^2+\cdots + c_{r+s}X^{r+s}$$ where $c_k=a_0b_k+a_1b_{k-1}+\cdots a_kb_0$ for $0\leq k\leq r+s.$
      
      This multiplication produces immediate examples of nonsingular bilinear maps
      $$\varPhi_0:\R^{8r+8}\times \R^{8s+8}\to \R^{8r+8s+8}$$ One simply identifies $\R^{8r+8}$ first with $\K^{r+1},$ and then visualizes a typical vector $\overrightarrow{a}=(a_0,a_1,\cdots, a_r )$ of $\K^{r+1}$ to be the polynomial $p(X),$ similarly with $\R^{8s+8}$ and $q(X).$ Then $\varPhi_0(p(X),q(X))$ is none other than $p(X)q(X).$ The nonsingularity of $\varPhi_0$ is tantamount  to the claim that $\K[X],$ like $\K,$ has no zero divisors. This can be easily proved by induction on the total degree $r+s.$ 
      
      It would be convenient to encode polynomial multiplication using a matrix scheme $M_0$ as follows

      \begin{equation}\label{M0}
      \xymatrix{  &  & & & &  & & & &M_0\\
      	& a_0b_0  \ar@{^-}[ur]  & a_0b_1  \ar@{^-}[ur]&     a_0b_2 \ar@{^-}[ur] & a_0b_3 \ar@{^-}[ur]& a_0b_4   \ar@{^-}[ur]    & \bullet \ar@{^-}[ur]& \cdots&\\
      L_0 \ar@{^-}[ur]	& a_1b_0 \ar@{^-}[ur] & a_1b_1  \ar@{^-}[ur] & a_1b_2  \ar@{^-}[ur] & a_1b_3  \ar@{^-}[ur] &  \bullet \ar@{^-}[ur]   &\bullet \ar@{^-}[ur]&\cdots &\\
      L_1 \ar@{^-}[ur] & a_2b_0 \ar@{^-}[ur]	 &   a_2b_1 \ar@{^-}[ur]  & a_2b_2 \ar@{^-}[ur] & a_2b_3 \ar@{^-}[ur] & \bullet \ar@{^-}[ur] &\bullet\ar@{^-}[ur]&\cdots\\
      L_2 \ar@{^-}[ur]	& a_3b_0 \ar@{^-}[ur]& a_3b_1 \ar@{^-}[ur] & a_3b_2 \ar@{^-}[ur]& a_3b_3 \ar@{^-}[ur] &\bullet \ar@{^-}[ur] &\bullet \ar@{^-}[ur] &\cdots\\
     L_3 \ar@{^-}[ur] & a_4b_0 \ar@{^-}[ur] & a_4b_1 \ar@{^-}[ur] & a_4b_2 \ar@{^-}[ur] & a_4b_3 \ar@{^-}[ur] & \bullet \ar@{^-}[ur] & \bullet \ar@{^-}[ur]&\cdots\\ 
  L_4 \ar@{^-}[ur] & a_5b_0 \ar@{^-}[ur] & a_5b_1 \ar@{^-}[ur] & a_5b_2 \ar@{^-}[ur] & a_5b_3 \ar@{^-}[ur] & \bullet \ar@{^-}[ur]  &\bullet\ar@{^-}[ur] &\cdots &\\
L_5 \ar@{^-}[ur]& L_6 \ar@{^-}[ur] & L_7 \ar@{^-}[ur] & L_8 \ar@{^-}[ur] &  \cdots \ar@{^-}[ur] & \bullet \ar@{^-}[ur] & \bullet \ar@{^-}[ur]}
      \end{equation}

      Here we put $a_{i-1}b_{j-1}$ at the $(i,j)^{th}$ entry of the matrix $M_0.$ The segmented lines $L_0,L_1,L_2$ etc. are successive  lines of slope 1, with unital decreases  in $y-$intercept, passing through various lattice points at which the $a_{i-1}b_{j-1}$ terms are situated. For instance, $L_2$ passes through $a_2b_0, a_1b_1$ and $a_0b_2.$ We shall use $\sum L_2$ to mean $a_2b_0+a_1b_1+a_0b_2$ for convenience. Indeed, the purpose of this encoding is to allow one to use $\sum L_k$ to symbolically stand for the coefficient $c_k$ of $X^k$ in the product $p(X)q(X)$ given above. When degree $p(X)= r,$ degree $q(X)= s,$ $M_0$ is of size $(r+1)\times (s+1),$ and can be used simultaneously to encode the bilinear map $\varPhi_0$ above as 
      $$\varPhi_0=\varPhi_{M_0}:\K^{r+1}\times \K^{s+1}\to \K^{r+s+1}$$ in the format
      $$\varPhi_{M_0}((a_0,a_1,\cdots, a_r),(b_0,b_1, \cdots, b_s))=(\sum L_0,\sum L_1,\cdots, \sum L_{r+s}).$$

\begin{remark}\label{remMbil}
	As a matter of fact, given any matrix $M$ of size $(r+1)\times (s+1), $ possibly $r=\infty$ or $s=\infty,$ of which each entry is an arbitrary bilinear form in $\overrightarrow{a}, \overrightarrow{b},$ a bilinear map $\varPhi_M$ will be automatically defined, with components  $\sum L_0,\sum L_1, \sum L_2, \cdots $ etc. The main effort of this paper is to seek out some $M's$ for which $\varPhi_M$ would be nonsingular.

There is no short supply of such $M's.$ We note that in the previous paragraph the guts of the induction argument for $\varPhi_0$'s nonsingularity is that, if $p(X)q(X)$ had vanishing leading coefficient $a_rb_s,$ then either $a_r=0$ or $b_s=0.$ The same conclusion, of course, also follows from $a_r\overline{b_s}=0,$ or from $-\overline{b_s}\overline{a_r}=0,$ and so on. This motivates us to bring in the following

\end{remark}

\begin{defi}
	 By a modification of the scheme $M_0,$ we mean a matrix $M$ of the same size, obtained from $M_0$ by replacing each $a_ib_j$ with any of the following choices 
	$$ \pm a_ib_j, \pm a_i\overline{b_j},\pm \overline{a_i}b_j, \pm \overline{a_i}\overline{b_j},$$
	$$\pm b_ja_i, \pm \overline{b_j}a_i, \pm b_j\overline{a_i},\pm \overline{b_j}\overline{a_i}.$$
\end{defi}
There is a total of 16 possibilities at each entry. Any modification $M$ of $M_0$ leads to a twisted polynomial multiplication $\odot_M$ for $\K[X]$ different from the traditional one, namely 

\begin{defi}
	The $M-$twisted product, or simply $M-$product, of $p(X)$ and $q(X)$ is defined to be 
	
	$$p(X)\odot_M q(X)=\hat{c}_0+\hat{c}_1X+\hat{c}_2X^2+\cdots$$ 
	where $\hat{c}_k=\sum L_k$ is the sum of all entries of $M$ falling on the line $L_k$ depicted as in the scheme $(\ref{M0}).$ 
\end{defi} 

The space of polynomials with coefficients in $\K$ under twisted product $\odot_M$ becomes a ring $\K_M[X]$ with not many nice properties. For example the constant polynomial $1$ may not be a 
two-sided multiplicative identity. However $\K_M[X]$ is still free of zero divisors, just like $\K[X](=\K_{M_0}[X]).$ An induction proof for the former works, almost verbatim, as it does for the latter. Thus we have 

\begin{theorem}\label{ThmM0modif}
	Let $M$ be any $(r+1)\times (s+1)$ matrix obtained from the matrix $M_0$ of equal size through modification. Then the bilinear map 
	$$\varPhi_M:\K^{r+1}\times \K^{s+1}\to \K^{r+s+1}$$
	defined by $\varPhi_M(p(X),q(X))=p(X)\odot_M q(X),$ or equivalently by 
	$$\varPhi_M(\overrightarrow{a},\overrightarrow{b})=(\sum L_0,\sum L_1, \cdots, \sum L_{r+s})$$ is again nonsingular.
\end{theorem}

\begin{remark}
	At first sight this theorem is not useful, as $\varPhi_M$ is of real type $(8r+8,8s+8,8r+8s+8),$ exactly the same as the type of $\varPhi_0 $ (or $\varPhi_{M_0}$). What, then, is the point of modifying?

But some modifications do indeed lead to novelty. The bilinearity of $\varPhi_{M_0}$ induces, in the most obvious way, an adjoint   	map 
$$adjM_0:\K^{r+1}-\{\overrightarrow{0}\}\to Mono_{\R}(\K^{s+1},\K^{r+s+1}),$$ where $Mono_{\R}$ means the space of monomorphisms from one real vector space to another. When $M_0$ is modified into $M,$ there is no reason why $adjM$ should be homotopic to $adjM_0.$ It is this  potentially new homotopy feature of $\varPhi_M$ that could be perhaps  exploited to produce new families of nonsingular bilinear maps. This strategy will be carried out in a case of $4\times 4$ matrices in Section 6. 
\end{remark}

\section{Neat factorization of polynomials and quasi irreducibility}\label{factorpol}

Let $M_0$ be the $(r+1)\times(s+1)$ matrix of Section \ref{twistedpolmult} that encodes standard multiplication in $\K[X]$ of one polynomial $p(X)$ of degree $\leq r$ with another polynomial $q(X)\leq s.$ Let $M$ be a modification of $M_0$ encoding a twisted multiplication $\odot_M$ of the same two polynomials. 

\begin{defi}\label{neatlyfactor}
	A polynomial $g(X)$ of degree $\leq r+s$ is said to be neatly factorized into an $M-$product of $p(X)$ and $q(X) $ if 
	\begin{enumerate}
		\item $g(X)=p(X)\odot_Mq(X),$
		\item  $p(X)$ and $q(X)$ are of positive degree and 
		\item The constant terms of $p(X)$ and $q(X)$ commute, i.e., $p(0)q(0)$ equals $q(0)p(0)$ as octonions. 
	\end{enumerate}
\end{defi}

\begin{defi}\label{quasiirr}
	A polynomial $g(X)$ is said to be quasi $M-$irreducible, if it does not admit any neat factorization in $\K_M[X].$ 
\end{defi}

As a matter of generality, this definition could be understood in the following way. One allows the $M_0$ in Section \ref{twistedpolmult} to have countable number of rows and columns, so that $M $ is also allowed  to be such; but one imposes the requirement that all modification of entries occur within the upper $(r+1)\times (s+1)$ block of $M_0,$ for some finite $r$ and finite $s.$ Throughout this paper, whether such generality ought to be in effect shall be clear from the context. 

Two mini examples will serve to illustrate    Definitions \ref{neatlyfactor} and \ref{quasiirr}.

\begin{example}\label{quadR}
	In $\R[X]$ usual multiplication of two linear polynomials into a quadratic is encoded by $$M_0=\begin{bmatrix}
	a_0b_0 & a_0b_1\\
	a_1b_0 & a_1b_1 
	\end{bmatrix}$$
	With respect to $M_0,$ $1+X^2$ is irreducible. If $M_0$ is modified into 
	
	$$M^-=\begin{bmatrix}
	a_0b_0 & a_0b_1\\
	a_1b_0 & -a_1b_1
	\end{bmatrix} $$
	then $1+X^2,$ as a quadratic in $\R_{M^-}[X],$ is factorizable, and of course any $M^-$ factorization is neat.

\end{example}

   \begin{example}\label{quadK}
   	Replace $\R$ by $\K$ in example \ref{quadR}, so that now all matrix entries are octonions. Let the $M_0 $ in \ref{quadR} be modified into $M$ where
   	
   	$$M=\begin{bmatrix}
   	a_0b_0 & a_0\ba   \\
   	\bz  a_1 & a_1b_1 
   	\end{bmatrix}$$  
   	
   	Let $i,j,k\in \H\subset \K$ be the usual imaginary quaternion units, considered here as octonions. Then one has an $\odot_M$ factorization of $i+iX^2$ into two linear factors, namely 
   	$$ i+iX^2=(j+iX)\odot_M (k+X).$$
   	Here  $a_0=j,a_1=i,b_0=k,b_1=1;\sum L_0=i,\sum L_1=0$ and $\sum L_2=i.$ Accordingly $i+iX^2$ is $M-$ reducible. 
   	But the above is not a neat factorization at all because the constant terms of the two factors do not commute: $jk\neq kj.$ In fact it is not hard to check here that $i+iX^2$ admits no neat factorization and is $M-$quasi irreducible. In hindsight, such quasi-irreducibility   is the key feature that facilitates the construction of nonsingular bilinear maps in \cite{Kee}. As we shall see, it also facilitates the new construction in Theorem \ref{Thmmain} below.

   \end{example}

\section{Twisted multiplication of two cubics}\label{twiscubics}

Traditional multiplication of two polynomials of at most cubic degree is encoded, as in Section \ref{twistedpolmult}, by the $4\times 4$ matrix $M_0=[a_ib_j], 0\leq i,j \leq 3,$ with octonion entries. In this section we consider a twisted multiplication $\odot_M$ given by a specific modification $M$ of $M_0,$ where

$$M=\begin{bmatrix}
a_0b_0 & a_0\ba    & a_0\bb    & \az b_3\\
\bz  a_1 & a_1b_1 & -\bb   a_1 & -b_3\aa \\
\bz  a_2 & a_2\ba    & a_2b_2 & a_2b_3\\
b_0a_3 & a_3b_1 & a_3b_2 & a_3b_3
\end{bmatrix}$$ 

\begin{theorem}\label{ThmCubic}
	For any nonzero $c\in \K$ and any strictly positive real number $\lambda,$ the polynomial $g(X)=c+\lambda cX^4$ in $\K_M[X]$ is quasi-irreducible. In other words $g(X)$ does not admit any neat $\odot_M$ factorization. 
\end{theorem}
\begin{proof}
	Since $g(X)$ has degree 4, its possible factorization must occur as either
	 
	Case 1: a linear times a cubic, or
	
	Case 2: a cubic times a linear, or 
	
	Case 3: a quadratic times another quadratic.
	
	Corresponding to these cases are the submatrices $M_1,M_2$ and $M_3$ of $M,$ of sizes $2\times 4, 4\times 2$ and $3\times 3$ respectively. $M_1$ is formed by $M's$ first two rows, $M_2$ by its first two columns, and $M_3$ is $M's$ principal $3\times 3$ sub-block. We display each one explicitly below, for easy tracking later.   
	 $$ M_1=\begin{bmatrix}
	 a_0b_0 & a_0\ba    & a_0\bb    & \az b_3\\
	 \bz  a_1 & a_1b_1 & -\bb   a_1 & -b_3\aa 
	 \end{bmatrix}$$
	
	$$M_2=\begin{bmatrix}
	a_0b_0 & a_0\ba    \\
	\bz  a_1 & a_1b_1 \\
	\bz  a_2 & a_2\ba    \\
	b_0a_3 & a_3b_1
	\end{bmatrix}\:\;\;\;\;\;\;
	M_3=\begin{bmatrix}
	a_0b_0 & a_0\ba    & a_0\bb    \\
	\bz  a_1 & a_1b_1 & -\bb   a_1 \\
	\bz  a_2 & a_2\ba    & a_2b_2
	\end{bmatrix}$$
	It suffices to establish that $g(X)$ is quasi $M_i-$irreducible for $i=1,2$ and $3.$ All cases are done by reduction to absurdity; that is, supposing $g(X)=p(X)\odot_{M_i}q(X)$ neatly and deriving therefrom a contradiction, $i=1,2,3.$ 
	
	Case 1. Suppose that, neatly, 
	$$c+\lambda cX^4=(a_0+a_1X)\odot_{M_1}(b_0+b_1X+b_2X^2+b_3X^3)$$
	The neatness requirement is 
	\begin{equation}\label{eqncommut}
	a_0b_0(=c)=b_0a_0\;\;\;\;\;\;\;\;\;\;\;\;\;\;\;\;\;\;\;\;\;\;\;\;\mbox{ (neatness) }
	\end{equation}
	
	Quick comparison of coefficients, using the scheme of $L-$lines of Section \ref{twistedpolmult} for the displayed $M_1,$ gives 
	
	\begin{eqnarray}
	\overline{b_0}a_1+a_0\overline{b_1}=0\label{E51}\\
	a_1b_1+a_0\overline{b_2}=0\label{E52}\\
	-\overline{b_2}a_1+\overline{a_0}b_3=0\label{E53}\\
	-b_3\overline{a_1}=\lambda c=\lambda a_0b_0, \lambda>0\label{E54}
	\end{eqnarray}
	
	where $c\neq 0$ entails $a_0\neq 0\neq b_0.$ It further entails $b_3\neq 0\neq a_1$ via (\ref{E54}) and also $b_2\neq 0$ via (\ref{E53}). That, in turn, forces $b_1\neq 0\neq a_1$ via (\ref{E52}). With obvious meaning for the notation evaluate $[b_0(\ref{E51})]b_1$ to obtain 
\begin{equation}\label{E55}
[b_0(\bz a_1)]b_1+[b_0(a_0\ba )]b_1=0
\end{equation}
	  Because $a_0$ and $b_0$ commute, restricted associativity shows that the second term on left equals $[(b_0a_0)\ba ]b_1$ which in turn equals $(b_0a_0)|b_1|^2,$ because $\ba$ and $b_1$ commute. Applying similar arguments to the first term we reduce (\ref{E55}) to  
	  \begin{equation}\label{E56}
	  |b_0|^2a_1b_1+|b_1|^2a_0b_0=0
	  \end{equation}
	  
	  Substituting (\ref{E52}) into (\ref{E56}) gives 
	  \begin{equation}
	  \label{E57} -|b_0|^2a_0\bb +|b_1|^2a_0b_0=0
	  \end{equation}
	   which confirms   $b_1\neq 0$ again.
	  
	  Left cancelling the nonzero $a_0$ factor from (\ref{E57}) one confirms $\bb$ to be a real multiple of $b_0,$ and obtains the crucial fact that $\bb$ commutes with $a_0$ since $b_0$ does.
	  
	  Next, evaluate $[a_0(\ref{E53})]\aa$ and simplify the result using restricted associativity to obtain 
	  \begin{equation}\label{E58}
	  -|a_1|^2a_0\bb +|a_0|^2b_3\aa =0
	  \end{equation}  
	  
	  Equations (\ref{E54}), (\ref{E57}) and (\ref{E58}) can be reorganized  into three homogeneous linear relationship, whit real coefficients, amongst the three  octonions $a_0b_0,$ $a_0\bb$ and $b_3\aa.$ Since $a_0b_0(=c)\neq 0,$ the $3\times 3$ matrix that collectively summarizes all such relation must have zero determinant. This matrix is 
	  
	  $$\begin{bmatrix}
	  -\lambda & 0 & -1\\
	  |b_1|^2 & -|b_0|^2 & 0\\
	  0 & -|a_1|^2 & |a_0|^2
	  \end{bmatrix}$$ 
	  Its determinant is $\lambda |a_0|^2|b_0|^2+ |b_1|^2|a_1|^2,$ which is strictly positive. We have thus arrived at  a contradiction.
	  
	  Case 2: This is very similar to case 1 and we can afford to be brief. Suppose that, neatly,
	  $$c+\lambda cX^4=(a_0+a_1X+a_2X^2+a_3X^3)\odot_{M_2}(b_0+b_1X).$$
	  
	  Then 
	  
	  \begin{eqnarray}
	  a_0b_0(=c)=b_0a_0 \mbox{ \;\;\;\;\;\;\;\;\;\;(neatness) }\label{E59}\\
	  \bz a_1+a_0\ba=0 \;\;\;\;\;\;\;\;\;\; \;\;\;\;\;\;\;\;\;\;\;\;\;\;\;\;\;\label{E510}\\
	  \bz a_2+a_1b_1=0\;\;\;\;\;\;\;\;\;\; \;\;\;\;\;\;\;\;\;\;\;\;\;\;\;\;\;\label{E511}\\
	  b_0a_3+a_2\ba =0\;\;\;\;\;\;\;\;\;\; \;\;\;\;\;\;\;\;\;\;\;\;\;\;\;\;\;\label{E512}\\
	  a_3b_1(=\lambda c)=\lambda a_0b_0, \lambda>0\label{E513}\;\;\;\;\;\;\;\;\;\; \;\;\;\;\;\;\;\;\;\;\;\;\;\;\;\;\;
	  \end{eqnarray}
	  
	     We now get, using $[b_0(\ref{E510})]b_1$ followed by substitution  
	  
	\begin{eqnarray}
	|b_0|^2a_1b_1+|b_1|^2a_0b_0=0 \label{E514}\\
	-|b_0|^2\bz a_2 +|b_1|^2a_0b_0=0 \label{E515}
	\end{eqnarray}  
	
	Use  $b_0(\ref{E515})$ to recognize $a_2$ as a real multiple of $b_0(a_0b_0),$ commuting, therefore, with $a_0b_0$ and  $\bz.$ This crucial commutativity helps produce, via $[\bz (\ref{E512})]b_1,$ the linear relationship
	\begin{equation}\label{E516}
	|b_0|^2 a_3 b_1+|b_1|^2\bz a_2=0
	\end{equation}
			
			Together with (\ref{E513}) and (\ref{E515}) there are three such relationships amongst $a_0b_0, \bz a_2,$ and $a_3b_1.$ Again $a_0b_0\neq 0$ forces the coefficient matrix to have zero determinant. This matrix is 
			
			$$\begin{bmatrix}
			-\lambda & 0 & 1\\
			|b_1|^2 & -|b_0|^2 & 0\\
			0 & |b_1|^2 & |b_0|^2
			\end{bmatrix}$$ 
			
		with determinant $\lambda |b_0|^4+|b_1|^4$ strictly positive. Contradiction !
		
		Case 3. With the $3\times 3$ matrix $M_3$ defining a twisted multiplication of two quadratics, the argument proceeds similarly, only that the route towards contradiction is a bit more devious. Supposing a neat factorization 
		
		$$c+\lambda cX^4=(a_0+a_1X+a_2X^2)\odot_{M_3}(b_0+b_1X+b_2X^2)$$
		leads to 
		
		\begin{eqnarray}
		a_0b_0(=c)=b_0a_0\mbox{ \;\;\;\;\;\;\;\;\;\;(neatness)} \label{E517}\\
		\bz a_1+a_0\ba =0 \;\;\;\;\;\;\;\;\;\; \;\;\;\;\;\;\;\;\;\;\;\;\;\;\;\;\;\label{E518}\\
		\bz a_2+a_1b_1+a_0\bb =0 \;\;\;\;\;\;\;\; \;\;\;\;\;\;\;\;\;\;\;\;\;\;\;\;\;\;\;\label{E519}\\
		a_2\ba -\bb a_1=0 \;\;\;\;\;\;\;\;\;\; \;\;\;\;\;\;\;\;\;\;\;\;\;\;\;\;\;\label{E520}\\
		a_2b_2(=\lambda c)=\lambda a_0b_0, \lambda >0 \label{E521}\;\;\;\;\;\;\;\;\;\; \;\;\;\;\;\;\;\;\;\;\;\;\;\;\;\;\;
		\end{eqnarray}
		where $c\neq0$ by hypothesis, thereby entailing  $a_2\neq 0 \neq b_2.$
		
		Exactly as before one can get 
		
		\begin{equation}
		\label{E522} |b_0|^2a_1b_1+|b_1|^2a_0b_0=0
		\end{equation}
		which shows $a_1b_1$ to be a real multiple of $a_0b_0,$ commuting, therefore, with $a_0$ and with $b_0.$ Now evaluate $[b_0(\ref{E519})]b_2$ to get 
		
		\begin{equation}
		\label{E523} |b_0|^2a_2b_2+[b_0(a_1b_1)]b_2+|b_2|^2a_0b_0=0
		\end{equation} 
		
		In (\ref{E523})  the middle term equals $[(a_1b_1)b_0]b_2$ which in turn equals $[a_1b_1]b_0b_2$ by restricted associativity. 
		
		Using (\ref{E522}) this term becomes 
		
		$$-[|b_0|^{-2}|b_1|^2a_0b_0](b_0b_2)$$
		which has $a_0b_0$ as a left factor. Observe that the first term of (\ref{E523}) has $a_0b_0$ as left factor too, on account of (\ref{E521}). Left cancelling this common factor reduces (\ref{E521}) to 
		
		\begin{equation}
		\label{E524} \lambda |b_0|^2-|b_0|^{-2}|b_1|^2b_0b_2+|b_2|^2=0
		\end{equation}   
		 
		 Since $b_0\neq 0$ and $\lambda > 0$ by hypothesis, (\ref{E524}) forces $b_1\neq 0$ and shows $b_0b_2$ to be a positive real number. Thus $b_2$ is just $\bz$ up to a real multiple, in resemblance with the conclusion from the earlier (\ref{E57}).  By properties of $\K$ listed in Section \ref{AsoOct}, $b_2\in \R(b_0)\subset\R(a_0,b_0).$ Equation (\ref{E521}) then implies $a_2\in \R(a_0,b_0)$ so that $a_2$ and $b_2$ commute. This crucial commutativity allows one to simplify $[b_2(\ref{E520})]b_1$ into
		 
		 \begin{equation}
		 \label{E525} |b_1|^2a_2b_2-|b_2|^2a_1b_1=0
		 \end{equation}   
		 
		 The octonions $a_0b_0,a_1b_1$ and $a_2b_2,$ with $a_0b_0\neq0,$ are now subject to homogeneous linear relations (\ref{E521}),(\ref{E522}) and (\ref{E525}) with real coefficients. Again the relevant coefficient matrix must have zero determinant. When that matrix is written down, like what was done in cases 1 and 2, its determinant is evaluated to be 
		 
		 $$|b_1|^2(\lambda|b_0|^2+|b_2|^2).$$
		 
		 Recalling $b_1\neq 0 $ from (\ref{E524}), we have reached our final contradiction. This completely establishes Theorem \ref{ThmCubic}

		\end{proof}
	
\begin{remark}
	\label{remtothmcubic1} The polynomial $g(X)=c+\lambda cX^4$ may or may not have an $M-$factorization that isn't neat. When $c=i,\lambda=1,$ one has 
	
	$$i+iX^4=(j+iX^2)\odot_M (k+X^2),$$
	an analogue of mini example \ref{quadK} in Section \ref{factorpol}.
	
	When $c=1,\lambda =1, 1+X^4$ has no factorization in the twisted polynomial ring $\K_M[X]$ whatsoever, because any $M-$factorization has to be neat, and thus contrary to Theorem \ref{ThmCubic}. 
\end{remark}	

\begin{remark}
	\label{remtothmcubic2} The effect of ``twisting'', i.e., conjugating, negating and factor transposing, becomes apparent if we recall that in the untwisted $\K[X]$ the following factorization is well known:
	
	$$1+X^4=(1+\sqrt{2}X+X^2)(1-\sqrt{2}X+X^2).$$
	Theorem \ref{ThmCubic} thus brings out the subtlety of octonionic arithmetic.  
\end{remark}

\section{The nonsingular bilinear map $\varPhi_{\widetilde{M}}:\K^4\times \K^4\to \K^7$ and octonion commutators}\label{Final}

We continue to study the $4\times 4$ matrix $M$ of Section \ref{twiscubics}. By Theorem \ref{ThmM0modif} it already defines a nonsingular map 

$$\varPhi_M:\R^{32}\times \R^{32}\to \R^{56}.$$

As it turns out $M$ can actually be adjusted slightly to become an $\widetilde{M}$ that defines better maps. To do so introduce the $4\times 4$ matrix,

$$N=\begin{bmatrix}
b_0a_0 & 0 & 0 & 0\\
0 & 0 & 0 & b_0a_0\\
0 & 0 & b_0a_0 & 0\\
0 & b_0a_0 & 0 & 0
\end{bmatrix}$$

We take the liberty to think of $N$ as encoding a very esoteric multiplication $\odot_N$ of cubic polynomials, namely

 $$(a_0+a_1X+a_2X^2+a_3X^3)\odot_N(b_0+b_1X+b_2X^2+b_3X^3)=b_0a_0+3b_0a_0X^4.$$
 
Alter $M$ to $\widetilde{M}=M-N.$ Each entry of $\widetilde{M}$ is bilinear in $\overrightarrow{a},\;\overrightarrow{b},$ so by Remark \ref{remMbil} it defines a bilinear  

$$\varPhi_{\widetilde{M}}:\K^4\times \K^4\to \K^7$$
where the first component of $\varPhi_{\widetilde{M}}(\overrightarrow{a},\overrightarrow{b})\in \K^7$ is an octonion commutator $a_0b_0-b_0a_0.$ This is the map of the section title.

\begin{theorem}
	\label{Thmmain} The bilinear map $\varPhi_{\widetilde{M}}$ is nonsingular. Also, because an octonion commutator has no real part, the type of $\varPhi_{\widetilde{M}}$ should be more accurately exhibited as 
	
	$$\varPhi_{\widetilde{M}}:\R^{32}\times \R^{32}\to \R^{55}.$$ lowering the range of $\varPhi_{\widetilde{M}}$ from $\R^{56}$ to $\R^{55}.$  
  \end{theorem}  
 
\begin{proof} 
The proof of nonsingularity is not hard thanks to the preparation in Section \ref{twiscubics}. The strategy is to suppose 

\begin{equation}
\label{E61} \varPhi_{\widetilde{M}}(\overrightarrow{a},\overrightarrow{b})=\overrightarrow{0}
\end{equation} 
and deduce that 

\begin{equation}
\label{E62} \mbox{ either } \overrightarrow{a}=\overrightarrow{0} \mbox{ or } \overrightarrow{b}=\overrightarrow{0}
\end{equation}
As is easily seen via the scheme of $L-$lines in Section \ref{twistedpolmult}, $\varPhi_{\widetilde{M}}=\varPhi_M-\varPhi_N, $ so that (\ref{E61}) means 

\begin{equation}
\label{E63} \varPhi_{M}(\overrightarrow{a},\overrightarrow{b})=\varPhi_N(\overrightarrow{a},\overrightarrow{b})
\end{equation} 

In terms of twisted polynomial product this is 

\begin{equation}
\label{E64} (a_0+a_1X+a_2X^2+a_3X^3)\odot_M(b_0+b_1X+b_2X^2+b_3X^3)=b_0a_0+3b_0a_0X^4
\end{equation}
 
 Comparing constant terms yields $a_0b_0=b_0a_0.$ We denote this common value by $c$ so that the right hand side of (\ref{E64}) reads $c+3cX^4.$
 
 Comparing coefficients in degrees $6$ and $5$ yields
 
 $$a_3b_3=0,\;a_3b_2+a_2b_3=0.$$
 This can happen only in one of the three circumstances below (cf. Section \ref{twiscubics}).  
 
 Case $(1)^\sim$  $a_2=0=a_3$ so that the left factor in (\ref{E64}) is at most linear, or
 
 Case $(2)^\sim$ $b_2=b_3=0.$ So that the right factor in (\ref{E64}) is at most linear, or
 
 Case $(3)^\sim$ $a_3=0,\;b_3=0.$ So that both factors are at most quadratic.
 
 We first show how to reach the desired conclusion (\ref{E62}) for Case $(3)^\sim.$   If both $\odot_M $ factors in (\ref{E64}) have positive degree, then (\ref{E64}) shows $c+3cX^4$ to be neatly $\odot_M$ factorizable into two quadratics, contrary to case 3 of Theorem \ref{ThmCubic}. Therefore either left or right factor must be a constant polynomial. Say left factor$=a_0,$ with $a_1=a_2=0$ (in addition to the case specification $a_3=0$). Directly from the encoding scheme of Section \ref{twistedpolmult}, the left hand side of (\ref{E64}) now reduces to $$a_0b_0+a_0\ba X+a_0\bb X^2+\az b_3X^3.$$ Coefficient comparison with right hand side gives $$a_0b_0=b_0a_0,\;a_0\ba=0,\;a_0\bb=0,\;\az b_3=0,\;0=3b_0a_0$$ This implies either $\overrightarrow{b}=\overrightarrow{0}$ or $a_0=0$ (and thus $\overrightarrow{a}=\overrightarrow{0}$), which is the desired (\ref{E62}).
 
 For the possibility that right factor $=b_0,$ with $b_1=b_2=0$ (in addition to the case specification $b_3=0$), the argument to reach (\ref{E62}) is entirely parallel.
 
 Finally, handling Cases $(1)^\sim$ and $(2)^\sim$ through appeals to Cases        1 and 2 of Theorem \ref{ThmCubic} proceeds in exact analogy with the $(3)^\sim$ case, and needs no further comment. The nonsingularity of $\varPhi_{\widetilde{M}}$ is now fully established.
 
 \end{proof}  

The question of existence of a bilinear map having same type as $\varPhi_{\widetilde{M}},$ with an octonion commutator in one component, was posed by J. Adem in \cite{Adem2}. Theorem \ref{Thmmain} answer this question in the affirmative. 

A good number of nonsingular bilinear maps now follow, some new, some  previously recorded. All are obtained from $\varPhi_{\widetilde{M}}$ by restrictions of domain and range. The choice of domain/range to restrict to is guided by the properties of octonionic commutators $a_0b_0-b_0a_0.$ We refer to \cite{Adem} or \cite{Adem1}  for a full     account of possible choices. Just for example, one can restrict $\varPhi_{\widetilde{M}}$ to

$$(\C\oplus\K^3)\times (\C\oplus \K^3)\to \{\overrightarrow{0}\}\oplus\K^6,$$
to obtain a type $\R^{26}\times \R^{26} \to \R^{48}$ originally envisaged by the first author as an extension of Adem's $\R^{18}\times \R^{18}\to \R^{32}$ in \cite{Adem}.

For yet another example, take $V\subset\K$ to be the $3-$dimensional real subspace spanned by the imaginary quaternionic units $i.j,k,$ with $5-$dimensional orthogonal complement $V^\bot\subset \K.$ Then because the commutator map for $\K$ restricts to $V^\bot \times V^\bot \stackrel{[,]}{\to} V,$ one obtains 

$$(V^\bot\oplus\K^3)\times (V^\bot\oplus\K^3)\to V\oplus\K^6$$     
to be another legitimate restriction of $\varPhi_{\widetilde{M}},$ leading to a new type $\R^{29}\times \R^{29} \to \R^{51}$ which is the most interesting among all possible restrictions. It generalizes the map $\R^{13}\times \R^{13} \to \R^{19}$ in the second author's Ph.D. thesis.  It also gives an immersion of $\R P^{28}$ into $\R^{50}$ without any need to use Postnikov obstruction theory.

We use  the following table and Propositions to summarize  the many nonsingular bilinear maps $\varPhi$ that can result from restricting $\varPhi_{\widetilde{M}}.$ Many of these are new, superseding,  for example, the maps constructed by Adem in \cite[Proposition 4.3]{Adem2}.

\begin{equation}
\label{table} \begin{tabular}{|c|c|c|c|c|c|c|c|c|}
\hline
h & 32 & 32 & 31 & 29 & 27 & 26 & 26 & 25\\
\hline
k & 32 & 26 & 27 & 29 & 27 & 30 & 26 & 32\\
\hline
 m & 55 & 54 & 53 & 51 & 49 & 52 & 48 & 48\\
 \hline
\end{tabular}
\end{equation}  

\begin{proposition}
	\label{Propexist}For each triple $(h,k,m)$ tabulated, there exits a nonsingular bilinear $\varPhi:\R^h\times \R^k\to \R^m$ obtained through restricting the domain and range of $\varPhi_{\widetilde{M}}$ in a suitable way.  
\end{proposition}

\begin{proposition}
	\label{Propfurthe} Additionally, by further restricting $\varPhi's$  domain and range, one can obtain a second level of nonsingular bilinear maps $\varPhi^-,$ of types $\R^{h-8}\times \R^k\to \R^{m-8}$ as well as $\R^h\times \R^{k-8}\to \R^{m-8}.$ These types either match, or supersede, all nearby  types documented so far in the literature 
\end{proposition}

\begin{proposition}
	\label{Propmorefurder} Further restriction of domain and range of $\varPhi^-$ yields a third level of nonsingular bilinear maps 
	$$\varPhi^{--}:\R^{h-8}\times \R^{k-8}\to \R^{m-16}$$ 
	We again leave out the details for selecting such restrictions, but point out that these $\varPhi^{--}$ coincide precisely with Adem's eight maps constructed in \cite[Theorem 3.6]{Adem}.  
\end{proposition}

In this sense Propositions  \ref{Propexist} and \ref{Propfurthe} become direct expansion of Adem's Theorem. Ultimately, to return the subject to its debut, one could further restrict the domain and range of $\varPhi^{--},$ to produce a lowest level of nonsingular bilinear maps. These are essentially the ones in \cite{Kee}.

Even though quite a number of maps constructed in this paper are new, in the range $h\leq32, k\leq 32$ there are maps in existing literature which supersede ours. One notable example is Milgram's $\R^{32}\times \R^{32}\to \R^{54}$ in \cite{Milgram}, reformulated by Adem in \cite{Adem2}. Milgram \cite{Milgram} has no explicit use of octonion commutators.

One naturally wonders whether Theorem \ref{Thmmain} can have higher dimensional analogue. For example, is there a nonsingular bilinear 

$$\K^8\times \K^8\to \K^{15}$$   
with a commutator component? An examination of the pattern of Proof in Section \ref{Final} and \ref{twiscubics} shows that, to get an answer one needs to struggle through a jungle of octonion  arithmetic, or to have new ideas. We leave this as an invitation to interested readers.


\begin{thebibliography}{1}
		
		\bibitem{Adem}Adem, J., Some immersions associated with bilinear maps, Bol. Soc. Mat. Mexicana 13(1968), 95-104.
		
		\bibitem{Adem1}Adem, J., On nonsingular bilinear maps. The Steenrod Algebra and its Applications,
		Lecture Notes in Math. 168(1970), Springer, 11-24.
		
		\bibitem{Adem2}Adem, J., On nonsingular bilinear maps II, Bol. Soc. Mat. Mexicana 16(1971), 64-70.
		
	
		
		\bibitem{Ginsburg}Ginsburg, M.,  Some immersions of projective space in Euclidean space, Topology 2(1963), 69-71.
		
	 
		
		
		\bibitem{Hopf}
		Hopf, H.,  Ein topologischer Beitrag zur reellen Algebra. Comment Math. Helv. 13(1941), 219-239.
		
			\bibitem{Kitchloo-Wilson}  Kitchloo, N., Wilson, S. The second real Johnson-Wilson theory and non-immersions of $\R P^n,$ Homology, Homotopy and Applications 10(3) 2008, 223-268. 
		
		\bibitem{Kee}Lam, K.Y., Construction of nonsingular bilinear maps,
		Topology 6(1967), 423-426.
		
		\bibitem{Kee0}Lam, K.Y., Construction of some nonsingular bilinear maps, Bol. Soc. Mat. Mexicana (2)
		13(1968), 88-94.
		
		\bibitem{Kee1} Lam, K.Y. On bilinear and skew-linear maps that are nonsingular, Quart. J. Math. Oxford
		(2) 19(1968), 281-288.
		
		
		\bibitem{Kee3}Lam, K. Y., Nonsingular bilinear maps and stable homotopy classes of spheres, Math. Proc.
		Cambridge Philos. Soc. 82(1977), 419-425.
		
	
		
		
		\bibitem{Kee5}Lam, K.Y., {\it Borsuk-Ulam Type Theorems and Systems of Bilinear Equations}. Geometry from the Pacfic Rim, Walter de Gruyer, Berlin-New York, 1997, pp. 183-194.
		
		\bibitem{Milgram}Milgram, J., Immersing projective spaces, Ann. of Math. 85(1967), 473-482. 
		
		
		\bibitem{Rees} Rees, E. G., Linear spaces of real matrices of large rank, Proc. Royal Soc. Edinburgh Sect. A 126(1996), No.1, 147-151.
		
		\bibitem{Steer} Steer, B., On the Embedding of projective spaces in euclidean space, Proc. London Math. Soc.(3) 21(1970), 489-501.
		
		\bibitem{Stiefel} Stiefel E., Uber Richtungsfelder in den projektiven Raumen und einen Satz aus der reellen
		Algebra, Comment Math. Helv. 13 (1941), 201-208. 
	
		
		
	
	
	\end{thebibliography}
\end{document}